\documentclass[11pt,a4paper]{article}
\usepackage{amsmath,amssymb}

\usepackage{amsmath}

\usepackage{color}
\usepackage[colorlinks,linkcolor=blue,citecolor=red]{hyperref}

\DeclareMathSymbol{\bbbr}{\mathalpha}{AMSb}{"52}
\DeclareMathSymbol{\bbbc}{\mathalpha}{AMSb}{"52}

\newtheorem{theorem}{Theorem}

\newtheorem{definition}[theorem]{Definition}

\newtheorem{lemma}[theorem]{Lemma}

\begin{document}

\title{ Confocal conics and 4-webs of maximal rank}

\author{{\Large Sergey I. Agafonov}\\
\\
Department of Mathematics,\\
S\~ao Paulo State University-UNESP,\\ S\~ao Jos\'e do Rio Preto, Brazil\\
e-mail: {\tt sergey.agafonov@gmail.com} }
\date{}
\maketitle
\unitlength=1mm

\vspace{1cm}

\begin{abstract} Confocal conics form an orthogonal net. Supplementing this net with one of the following: 1) the net of Cartesian coordinate lines aligned along the principal axes of conics,\\
 2) the net of Apollonian pencils of circles whose foci  coincide with the foci of conics,\\
  3) the net of tangents to a conic of the confocal family, we get a planar 4-web. We show that each of  these 4-webs is of maximal rank and characterize confocal conics from the web theory viewpoint.
\bigskip


\bigskip

\end{abstract}



\section{Introduction}
Confocal conics and quadrics, being objects of mathematical study from antiquity till modern times, are apparently well known and understood in both their most elementary (see \cite{HCV-32,GSO-16} ) and the most profound aspects (see \cite{DR-11}). So much surprising is a novel characterization \cite{BSST-17} thereof as coordinate surfaces (or lines)  of curvilinear orthogonal coordinate system $s_1,s_2,...s_n$ such that the Cartesian coordinates $x_1,x_2,...x_n$, aligned along the principal axes of quadrics/conics, are factored into  $n$ functions each depending only on one variable:
$$
x_i=f_{i,1}(s_1)f_{i,2}(s_2)...f_{i,n}(s_n),\ i=1,...,n.
$$
Considering foliations of $\mathbb{R}^n$ by Cartesian coordinate hyperplanes $x_i=const$ and curvilinear coordinate surfaces $s_k=const$, we obtain a (singular) $2n$-web. Taking logarithm of both sides of the above factorization formulae, we get $n$ Abelian relations of this $2n$-web. Recall that an Abelian relation of a $d$-web is a functional identity  of the simplest possible form
$$
\sum_{k=1}^du_k=0
$$
among some first integrals $u_i$ of the web foliations.

For $n=2$, we have 2 Abelian relations for the planar 4-web of confocal conics and Cartesian coordinate lines. It turns out that there is a third Abelian relation, independent of the above two and therefore the rank of the web (the number of non-trivial independent Abelian relations) is 3. This is maximal possible for a planar 4-web. It is well known that any planar 4-web of maximal rank is locally diffeomorphic to one formed by tangents to some algebraic curve of 4th class $C$. Let us call its dual $C^*$, which is a curve of degree 4, {\it the rank curve} or {\it rank quartic} of the 4-web of maximal rank. The rank curve can be reducible. If it has rectilinear components then the corresponding foliations of the 4-web linearization are pencils of lines. We describe confocal conics in terms of web theory as follows.
\begin{theorem}\label{CartesianTH}
A planar orthogonal net and coordinate lines of a Cartesian coordinate system  are a confocal family of conics and coordinate lines of the Cartesian coordinate system aligned along the principal axes of conics if and only if the following two conditions are satisfied:
\begin{enumerate}
\item  the 4-web formed by the net and coordinate lines of Cartesian coordinate system  is of maximal rank and
\item the rank quartic of the 4-web splits into a conic and 2 lines meeting at this conic, where the lines  correspond to the coordinate foliations of the  Cartesian coordinate system.
\end{enumerate}
\end{theorem}
Even though a 4-web of maximal rank is algebraizable (and therefore linerizable), to construct webs of maximal rank with prescribed geometry of its foliation leaves is always a challenge.
Another starting point for our study of confocal conics from the web-theoretical viewpoint is the following result \cite{Aw-18}, which is already formulated in the language of webs. Complete the net of confocal conics to a 4-web by 2 pencils of {\it Apollonian circles}: the {\it elliptic} pencil, formed by circles passing through the foci of confocal family, and the {\it hyperbolic} pencil, which is formed by circles orthogonal to all the circles from the elliptic pencil. Then any 3-subweb of the obtained 4-web is hexagonal \cite{Aw-18}. By the Theorem of Mayrhofer  and Reidemeister (see \cite{BB-38}), this 4-web is linearizable, all the foliations of its linearization being pencils of lines.
Observe that the elliptic and hyperbolic pencils can be viewed as coordinate lines of the so-called {\it bipolar} coordinate system. Then we describe confocal conics in a way similar to Theorem \ref{CartesianTH}  as  follows.
\begin{theorem}\label{BipolarTH}
For a 4-web of confocal family of conics and coordinate lines of bipolar coordinate system whose foci  coincide with foci of conics,  the following conditions hold true:
\begin{enumerate}
\item  the web rank is maximal,
\item the rank quartic splits into 4 lines in general position.
\end{enumerate}
Conversely, if a 4-web composed of an orthogonal net and of coordinate lines of bipolar coordinate system satisfies the above conditions and has only isolated singularities then some map of the form
$$
x+iy=z \mapsto  \frac{k\left(\frac{z+1}{z-1}\right)+1}{k\left(\frac{z+1}{z-1}\right)-1}=\bar{x}+i\bar{y},\ \ \ k \in \mathbb C^*,
$$
 brings it to a 4-web of confocal family of conics and of coordinate lines of bipolar coordinate system whose foci coincide with the foci of conics.
\end{theorem}
Note that all the orthogonal coordinate systems considered so far are locally conformally equivalent to a Cartesian one. For brevity, we will call the corresponding net of coordinate lines {\it coformally flat}.

Finally, we consider a 4-web formed by confocal conics and tangent lines to one conic from the confocal family. Though it seems to be never formulated explicitly, analytic treatment of the net of tangents implies immediately (see for instance \cite{BSST-17} or \cite{IT-17}) that this 4-web is parallelizable, in particular, the two 3-subwebs, whose two foliations are tangents, are hexagonal. We show that, in nondegenerate cases, the  last condition characterizes the confocal family of conics.
Namely, take a {\it linear} net $\cal L$ composed of straight lines in a planar domain and define its {\it bisector net} as being formed by integral curves of directions bisecting the angles of line intersection. Then, for thus constructed 4-web, we prove the following result.

\begin{theorem}\label{TangentTH}
For the bisector net $\cal N$ of a linear net $\cal L$, the following conditions are equivalent:
\begin{enumerate}
\item the net $\cal L$ form a  hexagonal 3-web with one of the foliations of  $\cal N$,
\item the net $\cal N$ is conformally flat.
\end{enumerate}
Moreover, under the above two restrictions
\begin{itemize}
\item the rank of the 4-web formed by $\cal N$ and $\cal L$  is maximal,
\item the rank quartic of this 4-web splits into 4 concurrent lines forming a harmonic quadruplet,
\item the lines of ${\cal L}$ are either tangent to some smooth conic $Q$ or form 2 pencils,
\item the leaves of the net $\cal N$ are arcs of  some  confocal family of conics that includes  $Q$, if $Q$ is not a circle,
\item if $Q$ is a circle then the net $\cal N$ is formed by coordinate lines of a polar coordinate system, $Q$ being one of these lines.
\end{itemize}

\end{theorem}

\section{Web rank and linearization of 4-webs of maximal rank}

In this section, we recall concepts and facts necessary for our study.  The reader may turn to the classical monograph \cite{BB-38} or to the modern treatment \cite{PP-15} for more detail. All the considerations are local, webs and their foliations are defined on some planar open set $U$.

We will capture the foliations ${\cal F}_i$ of our planar webs as integral curves of ordinary differential equations (ODE) $\omega_i=0$, where $\omega_i$ are differential 1-forms. Clearly, these forms are defined up to multiplication by non-vanishing functions.
\begin{definition} A $d$-tuple of differential 1-forms $(\omega_1,...,\omega_d)$ is an Abelian relation of a $d$-web ${\cal W}_d$ with foliations ${\cal F}_i$, $i=1,2,...,d$ if
\begin{enumerate}
\item the leaves of ${\cal F}_i$ are solutions of ODEs $\omega_i=0$  (the choice $\omega_i\equiv 0$ is allowed),
\item the forms $\omega_i$ are closed: $d\omega_i=0$,
\item holds $\sum_{i=1}^d\omega_i=0$.
\end{enumerate}
The Abelian relations of the $d$-web ${\cal W}_d$ form a vector space. Its dimension is called the rank  of the web ${\cal W}_d$.
\end{definition}
Bol showed \cite{Bn-32} that the rank is finite and obtained the sharp bound $rk({\cal W}_d)\le \frac{1}{2}(d-1)(d-2)$. In particular, a 3-web is hexagonal if and only if its rank is maximal and is equal to one, i.e. there is a non-trivial Abelian relation. The rank is maximal for the so-called {\it algebraic webs}, obtained as a linear web of tangents to some algebraic curve $C$ of class $d$ in the plane.  (Recall that the class of an algebraic curve $C$ is the degree of its dual $C^*$. Thus, the term "Abelian" comes from the celebrated Abel theorem, which, in fact, implies the existence of the space of Abelian relations with desired dimension \cite{BB-38,CG-78}.) The planar 4-webs of maximal rank stand out by 2 facts: 1) they are algebraizable, i.e. locally diffeomorphic to algebraic ones (this is not true for $d\ge 5$) and 2) their linearizations (i.e.  local diffeomorphisms rectifying the leaves of all web foliations) are unique up to  projective transformations (this is not true for 3-webs).
Thus, for any planar 4-web of maximal rank ${\cal W}_4$, there is a linearization mapping the web leaves to tangents   to a curve dual to some quartic. For brevity, we will call the projective class of this quartic the {\it rank curve} or {\it rank quartic} of ${\cal W}_4$. The first remarkable fact above is due to the Sophus Lie studies of surface of double translations \cite{Lf-82}. In fact, he also provided a construction for a representative of the rank quartic. Essentially, it is as follows: for each foliation ${\cal F}_i$ choose a first integral $s_i$ and write a basis of Abelian relations as $\sum_{i=1}^4x_i^{\alpha}(s_i)ds_i=0$, $\alpha =1,2,3$ (recall that the rank of our web is 3), then the four points $[x^1_i:x^2_i:x^3_i]\in \mathbb P^2$ trace four arcs $C_i^*$ of the rank quartic (see, for example, \cite{BB-38} or \cite{Cw-82}). The quartic $C^*$ is allowed to be reducible. For instance, if it has a straight line as a component then the corresponding foliation is a pencil of lines, i.e. the corresponding component of the curve $C$ degenerates into a point.

\section{Confocal conics and Cartesian coordinates}

\noindent {\bf Proof of Theorem \ref{CartesianTH}:} Leaves of an orthogonal net can be described as integral curves of the ODEs $\omega_i=0$ with
\begin{equation}\label{forms}
\omega_1=Td(x)+d(y),\ \ \ \ \omega_2=-d(x)+Td(y).
\end{equation} Differentiating and excluding parameters in the family of confocal conics
\begin{equation}\label{confocal}
\frac{x^2}{a^2-\lambda}+\frac{y^2}{b^2-\lambda}=1,
\end{equation}
one gets the following two equations for the slope $T$:
\begin{equation}\label{Tconfocal}
T_x= \frac{T(y+2xT-yT^2)}{xy(T^2+1)}, \ \ \ \ T_y=\frac{ T(x-2yT -xT^2)}{xy(T^2+1)}.
\end{equation}
Let us look for Abelian relations in the form
\begin{equation}\label{Abel}
-Mdx-Ndy+K\omega_1+L\omega_2=0.
\end{equation}
Then it is immediate that
\begin{equation} \label{MN}
M=KT-L, \ \  \ \ \  N=LT+K.
\end{equation}
 Among the four conditions
\begin{equation} \label{closed}
M_y=N_x=0,\ \ \ \ d(K\omega_1)=d(L\omega_2)=0
\end{equation}
only three are independent. Introduce $J$ by $2J=K_y-L_x$ and find $K_x,K_y,L_x,L_y$ as functions of $x,y,T,K,L,J$:
\begin{equation}
\begin{array}{l}
K_x=L_y= TJ-\frac{T(xT^2+2yT-x)}{2xy(T^2+1)}K+\frac{T(yT^2-2xT-y)}{2xy(T^2+1)}L,\\
\\
K_y= J+\frac{(xT^2+2yT-x)}{2xy(T^2+1)}K+\frac{(yT^2-2xT-y)}{2xy(T^2+1)}L,\\
\\
L_x= K_y-2J.\\
\end{array}
\end{equation}
The compatibility conditions $d(K_xdx+K_ydy)=d(L_xdx+L_ydy)=0$ give $J_x$ and $J_y$:
\begin{equation}
\begin{array}{l}
J_x= \frac{P_1}{2xy(T^2+1)^2}J-\frac{P_2}{4x^2y^2(T^2+1)^4}K+\frac{P_3}{4x^2y^2(T^2+1)^4}L,\\
\\
J_y= -\frac{Q_1}{2xy(T^2+1)^2}J+\frac{Q_2}{4x^2y^2(T^2+1)^4}K-\frac{Q_3}{4x^2y^2(T^2+1)^4}L,
\end{array}
\end{equation}
where the polynomials $P_i,Q_i$ are as follows:
\begin{equation}
\begin{array}{l}{\scriptstyle
P_1=xT^5+7yT^4-14xT^3-10yT^2+xT-y,\ \ \ \ Q_1=yT^5-7xT^4-14yT^3+10xT^2+yT+x,}\\
\\
{\scriptstyle P_2=x^2T^9+9xyT^8+14(y^2-x^2)T^7-20xyT^6+(10y^2-16x^2)T^5-66xyT^4+(30x^2-38y^2)T^3+}\\
\\
{\scriptstyle \ \ \ \ +28xyT^2-(x^2+2y^2)T+xy},\\
\\
{\scriptstyle Q_2=xyT^9+(2y^2-7x^2)T^8-16xyT^7-(2x^2+10y^2)T^6-38xyT^5+(36x^2-42y^2)T^4+40xyT^3+}\\
\\
{\scriptstyle \ \ \ \ +2(x^2+y^2)T^2-3xyT+3x^2},\\
\\
{\scriptstyle P_3=xyT^9+(7y^2-2x^2)T^8-16xyT^7+(10x^2+2y^2)T^6-38xyT^5+
(42x^2-36y^2)T^4+40xyT^3+}\\
\\
{\scriptstyle \ \ \ \ -2(x^2+y^2)T^2-3xyT-3y^2},\\
\\
{\scriptstyle Q_3=y^2T^9-9xyT^8+14(x^2-y^2)T^7+20xyT^6+(10x^2-16y^2)T^5+66xyT^4xy+(30y^2-38x^2)T^3+}\\
\\
{\scriptstyle \ \ \ \ -28xyT^2-(2x^2+y^2)T-xy.}\\
\end{array}
\end{equation}
Now the form $J_xdx+J_ydy$ with the above expressions for $J_x$ and $J_y$ is closed and the rank of our web is 3, which is maximal possible for planar 4-webs.

Let us choose some basis for Abelian relations. The corresponding 3 solutions $(J_{\alpha},K_{\alpha},L_{\alpha})$,  $\alpha=1,2,3$ of  the above derived integrable Frobenius system define 3 maps
 $\vec{J},\vec{K}, \vec{L}: U\to \mathbb R^3$, where
 $\vec{J}=(J_1,J_2,J_3)^t$, $\vec{K}=(L_1,L_2,L_3)^t$, $\vec{L}=(L_1,L_2,L_3)^t$. For all $p\in U$, where $U$ is the common domain of the solutions, their images $\vec{J}(p),\vec{K}(p), \vec{L}(p)$ constitute a basis, any vector $\vec{v} \in \mathbb R^3$ can be represented as
 $\vec{v}=X\vec{K}+Y\vec{L}+Z\vec{J}$. With $\vec{v}$ fixed, the coordinates $X,Y,Z$ are functions of $p\in U$. One easily computes their differentials from the equation
 $$
 d\vec{v}=dX\vec{K}+dY\vec{L}+dZ\vec{J}+Xd\vec{K}+Yd\vec{L}+Zd\vec{J}=0,
 $$
applying the found expressions for the derivatives of $J,K,L$.

Now one checks that
$d(L_1)=0$ if $L_1=0$, where
$$
{\scriptstyle L_1(X,Y,Z)=X+TY-\frac{(yT^5-3xT^4-6yT^3+6xT^2+yT+x)}{2xy(T^2+1)^2}Z.}
$$
Considering $X,Y,Z$ as homogeneous coordinates on the plane $\mathbb{P}^2$, one concludes that the equation $L_1=0$ defines a straight line. A direct computation also shows that the points  $\vec{M}=T\vec{K}-\vec{L}$ and $\vec{M}_x$ lie on this line.
Thus, the rank curve arc, corresponding to the foliation $dx=0$ is a straight line.

Similarly, the rank curve arc, corresponding to the foliation $dy=0$ is also a straight line, whose equation reads as
$$
{\scriptstyle L_2(X,Y,Z)=-TX+Y-\frac{(xT^5+3yT^4-6xT^3-6yT^2+xT-y)}{2xy(T^2+1)^2}Z=0,}
$$
since $\vec{N}=\vec{K}+T\vec{L}$ and $\vec{N}_y$ lie on this  line.

Finally, the equation
$$
{\scriptstyle XY-\frac{(yT^2-2xT-y)}{2xy(T^2+1)}XZ+\frac{(xT^2+2yT-x)}{2xy(T^2+1)}YZ -  \frac{[xy(T^8-4T^6-26T^4+12T^2+1)+2(y^2-x^2)T(T^6+T^4-9T^2-1)]}{4x^2y^2(T^2+1)^4}Z^2=0}
$$
defines a smooth conic $c$ such that the curves parametrized by $[\vec{K}],[\vec{L}]$ (as $p$ ranges over $U$) have the second order contact to $c$ at any point. Therefore the rank curve arcs, corresponding to the remnant two foliations of the web, belong to the conic $c$. The intersection of the lines $L_1=0$ and $L_2=0$ is on this conic.
Thus, we have proved the necessary part of the theorem.

To prove the converse claim, observe that the web satisfying the conditions of the theorem is linearizable. Moreover, its rank curve is projectively equivalent to the one of a 4-web of confocal  conics and coordinate lines of Cartesian coordinate system aligned along the principal axes of conics.
Thus, the linearizations of these two webs are also projectively equivalent, in particular, the webs are locally diffeomorphic.
By this local diffeomorphism, the Cartesian coordinate lines map to Cartesian coordinate lines, hence  the diffeomorphism is of the form $(x,y)\to (u(x),v(y))$ (after interchanging $x$ and $y$, if necessary). Since it respects orthogonality of two pairs of directions at each point, it is conformal. Therefore $u(x)=kx+a$ and $v(y)=ky+b$ and the webs are related by homothety and translation.
\hfill $\Box$\\

The above approach permits  to describe more general orthogonal coordinate systems.
\begin{lemma}There exist isothermal coordinates such that the tangent vectors to their coordinate lines are kernels of the forms (\ref{forms}) if and only if the slope $T$ satisfies the equation
\begin{equation}\label{laplasT}
T_{xx}+T_{yy}=\frac{2T(T^2_x+T^2_y)}{1+T^2}.
\end{equation}
\end{lemma}
  {\it Proof:} The existence of the desired coordinates amounts to existence of a common integrating factor for the forms (\ref{forms}), which is equivalent to (\ref{laplasT}). \hfill $\Box$\\

\noindent{\bf Remark 1.} Equation (\ref{laplasT}) is the Laplace equation for the isothermal coordinates. For brevity, we will call an orthogonal net, satisfying the above Lemma, {\it conformally flat}.

\begin{theorem}\label{twoconf}
Suppose that two orthogonal nets form a 4-web of maximal rank in some planar open domain and one net is conformally flat then the other is also conformally flat and the rank curve splits into two conics.
\end{theorem}
{\it Proof:} Since one of the net is conformally flat, we can conformally map it into the coordinate net of some Cartesian coordinate system.  Consider again the forms (\ref{forms}), where $x,y$ are these Cartesian coordinates, and the Abelian relations (\ref{Abel}) (but without imposing equations (\ref{Tconfocal})). Again one gets $M,N$ as in (\ref{MN}). Then conditions (\ref{closed}) give
\begin{equation}
\begin{array}{l}
K_x=L_y=TJ+\frac{T_y}{2}K-\frac{Tx}{2}L,\\
\\
K_y= J-\frac{Ty}{2T}K-\frac{Tx}{2T}L,\ \ \ \ L_x=K_y-2J,
\end{array}
\end{equation}
and from $d(K_xdx+K_ydy)=d(L_xdx+L_ydy)=0$ we compute
\begin{equation}
\begin{array}{l}
J_x= -\frac{[(5T^2+1)T_x+(T^3-3T)T_y]}{2T(T^2+1)}J-\frac{ [(2T^3-2T)T_{xy}+(3T^2+3)T_xT_y+(T^3-3T)T^2_y]}{4T^2(T^2+1)}K+\\
\\
\ \ \  +\frac{[(2T^3+2T)T_{xx}-4T^2T_{xy}+(T^3+T)T_xT_y+(T^2-3)T^2_x]}{4T^2(T^2+1)}L,\\
\\
J_y= \frac{[(T^3-3T)T_x-(5T^2+1)T_y]}{2T(T^2+1)}J+\frac{[-(2T^3+2T)T_{yy}-4T^2T_{xy}+(T^3+T)T_xT_y-(T^2-3)T^2_y] }{4T^2(T^2+1)}K+\\
\\
\ \ \  +\frac{[(2T^3-2T)T_{xy}+(3T^2+3)T_xT_y-(T^3-3T)T^2_x]}{4T^2(T^2+1)}L.
\end{array}
\end{equation}
The web has rank 3 if and only if the coefficients of $J,K,L$ in  the equation $d(J_xdx+J_ydy)=0$ all vanish. The coefficient of $J$ yields (\ref{laplasT}) thus implying the conformal flatness.  Let us define $S$ by
$$
T_{yy}=-S+\frac{T(T_x^2+T_y^2)}{T^2+1},\ \ \ \ T_{xx}=S+\frac{T(T_x^2+T_y^2)}{T^2+1}.
$$
Then $d(J_xdx+J_ydy)=0$ gives $S_x$ and, consequently, $S_{xy}$. Further, $d(d(T_{xy}))=0$ determines $S_{yy}$ and, finally, $d(d(S_y))=0$ fixes $S_y$. Thus, we have
\begin{equation}
\begin{array}{l}\label{SxSy}
S_x=\frac{T_{xy}[4TT_x+(T^2-3)T_y]}{T(T^2+1)}+\frac{T_xT_y[(T^2+3)T_y-4TT_x]}{T^2(T^2+1)}+\frac{ST_x(3T^2+1)}{T(T^2+1)},\\
\\
S_y=\frac{T_{xy}[4TT_y-(T^2-3)T_x]}{T(T^2+1)}-\frac{T_xT_y[(T^2+3)T_x+4TT_y]}{T^2(T^2+1)}+\frac{ST_y(3T^2+1)}{T(T^2+1)}.
\end{array}
\end{equation}
With these $S_x$ and $S_y$, holds  $d(S_xdx+S_ydy)=0$. Therefore the exterior system for $K,L,J,$ $T,T_x,T_y,T_{xy},S$ is integrable and the web has the maximal rank 3.

To study the rank curve of this 4-web, we use again the moving frame $\vec{J},\vec{K}, \vec{L}$.  The arcs corresponding to the foliations $dx=0$ and $dy=0$ lie on the conic
\begin{equation}\label{conicdxdy}
\begin{array}{l}{\scriptstyle
X^2-Y^2+\frac{T^2-1}{T}XY+\frac{(T^4-4T^2-1)T_x+2T(1-T^2)T_y}{2T^2(T^2+1)}XZ-\frac{(T^4-4T^2-1)T_y-2T(1-T^2)T_x}{2T^2(T^2+1)}YZ+}\\
\\
\ \ \ {\scriptstyle -\frac{T(T^4-6T^2-3)(T_x^2-T_y^2)+(T^6-5T^4-T^2-3)T_xT_y-2T(T^4-1)T_{xy}+4ST^2(T^2+1)}{4T^3(T^2+1)^2}Z^2=0.}
\end{array}
\end{equation}
Therefore the rank curve, being a quartic, is a union of two conics. \hfill $\Box$\\

\noindent{\bf Remark 2.} In general, the above two conics of the rank curve are smooth.
The conic (\ref{conicdxdy}) splits into two lines if and only if
$$
S=\frac{T(T_x^2-T_y^2)}{T^2+1}-\frac{(T^4-4T^2-1)T_xT_y}{2T^2(T^2+1)}+\frac{(T^2-1)T_{xy}}{2T},
$$
the line equations being
\begin{equation}
\begin{array}{l}
{\scriptstyle TX-Y-\frac{(T^2+1)T_x+T(T^2-3)T_y}{2T(T^2+1)}Z=0,}\\
\\
{\scriptstyle X+TY+\frac{T(T^2-3)T_x-(T^2+1)T_y}{2T(T^2+1)}Z=0.}
\end{array}
\end{equation}
Moreover, the above expression for $S$ is compatible with (\ref{SxSy}). The lines meet in the other component of the rank quartic, which is the conic
\begin{equation}
{\scriptstyle XY+\frac{T_x}{2T}XZ-\frac{T_y}{2T}YZ -\frac{(T^2+3)T_xT_y-2TT_{xy}}{4T^2(T^2+1)}Z^2=0},
\end{equation}
if and only if
$$
T_{xy}=\frac{8T^2(T_y^2-T_x^2)}{(T^2+1)^3}+\frac{(T^6+11T^4-5T^2+1)T_xT_y}{T(T^2+1)^3}.
$$
This restriction for $T_{xy}$ also does not impose further constraints for the exterior system that is now  reduced to the above equation for $T_{xy}$ and the following two equations:
$$
T_{xx}=\frac{2T(T^4+3)T_x^2}{(T^2+1)^3} +\frac{2(3T^4-2T^2+3)T_xT_y}{(T^2+1)^3} +\frac{4T(T^2-1)T_y^2}{(T^2+1)^3},
$$
$$
T_{yy}= \frac{4T(T^2-1)T_x^2}{(T^2+1)^3}   -\frac{2(3T^4-2T^2+3)T_xT_y}{(T^2+1)^3} + \frac{2T(T^4+3)T_y^2}{(T^2+1)^3}.
$$
Thus we obtained an integrable Frobenius system for $T,T_x,T_y$ describing a 4-web whose rank curve is projectively equivalent to the 4-web of confocal family of conics and coordinate lines of Cartesian coordinate system aligned along the principal axes of conics.\\

\smallskip

 \noindent{\bf Remark 3.} The integral curves of the foliation $Tdx+dy=0$ are not conics for the 4-web, defined by the general solution to the Frobenius system for $T,T_x,T_y,T_{xy},S$ derived in the proof of Theorem \ref{twoconf}. One can check this as follows: the curves are solutions of the ODE $\frac{dy}{dx}=-T$ and they do not satisfy the ODE for general conic $\frac{d^3}{dx^3} \left(\frac{d^2y}{dx^2}\right)^{-2/3}=0$.

 One may capture the curves of the non-Cartesian net as level curves of the real and imaginary parts $u,v$ of some (possibly singular) conformal map $g(z)=g(x+iy)=u+iv$, sending the net into a Cartesian net $u=const,\ v=const$. This $g$ is a solution of an ODE.  One has $\frac{dg}{dz}=\partial_x(u+iv)=u_x+iv_x=RT-iR$, where $du=R(Tdx+dy)$ and $dv=R(-dx+Tdy)$ for some common integrating factor $R$. From the differential system  for $T$ and $S$ and the equation for $R$, which reads as
 $$
d\ln (R)=-\frac{1}{2}d\ln (1+T^2)+\frac{T_ydx-T_xdy}{1+T^2},
 $$
one can derive that ODE for $g$.

\section{Confocal conics and  bipolar coordinates}

As we have seen, the rank curve of a 4-web composed of two orthogonal conformally flat nets fixes this 4-web up to some conformal map preserving one of the nets. In this section, we need a conformal map preserving the coordinate foliations of a bipolar coordinate system. In the complex language, the  maps, preserving the coordinate foliations of the bipolar coordinate system with foci at $(\pm 1,0)$,  are exhausted by
\begin{equation}\label{preservebipolar}
x+iy=z \mapsto \frac{k\left(\frac{z+1}{z-1}\right)^{\alpha}+1}{k\left(\frac{z+1}{z-1}\right)^{\alpha}-1},\ \ \ k \in \mathbb C^*,\ \alpha \in \mathbb R^*.
\end{equation}
(To see this, one  applies a fractional-linear map to send $-1$ to $0$ and $1$ to $\infty$ and observes that the maps preserving the coordinate foliations of the polar coordinates are $z\mapsto kz^{\alpha}$.)
\medskip

\noindent {\bf Proof of Theorem \ref{BipolarTH}:} Let the forms (\ref{forms}) determine the
foliations by confocal conics. Then, as was shown above, holds (\ref{Tconfocal}). For the forms
\begin{equation}\label{Sforms}
\omega_1=Sd(x)+d(y),\ \ \ \ \omega_2=-d(x)+Sd(y),
\end{equation}
corresponding to the foliations by coordinate lines of the bipolar coordinate system whose foci coincide with the foci of conics, one gets, by excluding the parameter via differentiation, the following expression for  $S$:
$$
S=\frac{2T}{1-T^2}.
$$
Any 3-subweb of our 4-web is hexagonal (see \cite{Aw-18}). Then by the Theorem of Mayrhofer  and Reidemeister (see \cite{BB-38}, p.82), this 4-web is linearizable, all the foliations of its linearization being pencils of lines. Thus, the rank curve is a union of 4 lines. To clarify their mutual positions, we again calculate the compatibility conditions of equations for the Abelian relation (\ref{Abel}). These conditions read as
$$
M= \frac{(1-T^2)(K-TL)}{1+T^2},\ \ \ \ \ N= \frac{(1-T^2)(KT+L)}{1+T^2},
$$
\begin{equation}
\begin{array}{l}
K_x= TJ-\frac{T(xT^4+2yT^3+4xT^2+6yT-x)}{2xy(T^2+1)^2}K+\frac{T(yT^4-2xT^3-4yT^2+2xT-y)}{2xy(T^2+1)^2}L,\\
\\
K_y= J+\frac{(xT^4+2yT^3-4xT^2-2yT-x)}{2xy(T^2+1)^2}K+\frac{(yT^4-2xT^3-4yT^2+2xT-y)}{2xy(T^2+1)^2}L,\\
\\
L_x=K_y-2J,\\
\\
L_y= TJ-\frac{T(xT^4+2yT^3-4xT^2-2yT-x)}{2xy(T^2+1)^2}K+\frac{T(yT^4-2xT^3+4yT^2-6xT-y)}{2xy(T^2+1)^2}L,\\
\\
J_x=\frac{(xT^5+7yT^4-18xT^3-18yT^2+5xT-y)}{2xy(T^2+1)^2}J+\\
\\ \ \  -\frac{(x^2T^7+9xyT^6+(14y^2-15x^2)T^5-37xyT^4+(11x^2-16y^2)T^3+19xyT^2+(2y^2-5x^2)T+xy)}{4x^2y^2(T^2+1)^3}K+\\
\\
\ \ +\frac{(xyT^7+(7y^2-2x^2)T^6-25xyT^5+(24x^2-17y^2)T^4+35xyT^3+(5y^2-6x^2)T^2-3xyT-3y^2)}{4x^2y^2(T^2+1)^3}L,\\
\\
J_y=-\frac{(yT^5-7xT^4-18yT^3+18xT^2+5yT+x)}{2xy(T^2+1)^2}J+\\
\\
+\ \ \frac{(xyT^7+(2y^2-7x^2)T^6-25xyT^5+(17x^2-24y^2)T^4+35xyT^3+(6y^2-5x^2)T^2-3xyT+3x^2)}{4x^2y^2(T^2+1)^3}K+\\
\\
\ \ -\frac{(y^2T^7-9xyT^6+(14x^2-15y^2)T^5+37xyT^4+(11y^2-16x^2)T^3-19xyT^2+(2x^2-5y^2)T-xy)}{4x^2y^2(T^2+1)^3}L.
\end{array}
\end{equation}
Using the moving frame $\vec{J},\vec{K}, \vec{L}$, we calculate the lines of the rank curve.
For the foliations by circles, they are
$$
TX+Y+\frac{(xT+y)}{2xy}Z=0,\ \ \ \ X-TY+\frac{(yT-x)}{2xy}Z=0,
$$
and for the foliations by conics
$$
X+\frac{(xT^4+2yT^3-4xT^2-2yT-x)}{2xy(T^2+1)^2}Z=0, \ \ \ \ Y-\frac{(yT^4-2xT^3-4yT^2+2xT-y)}{2xy(T^2+1)^2}Z=0.
$$
They are in general position.

Now consider the 4-web ${\cal W}_b$ composed by coordinate lines of the bipolar coordinate system with foci at $(\pm 1,0)$ and by confocal conics with foci at these points and some 4-web  $\cal W$, satisfying the hypothesis of the converse part of Theorem \ref{BipolarTH} and having the same foci. Local linearizations of ${\cal W}_b$ and  $\cal W$ are projectively equivalent since their rank curves are projectively equivalent. Therefore there is a local diffeomorphism $\varphi$  mapping leaves of ${\cal W}_b$, intersecting a neighbourhood of some non-singular point $p_0$ of ${\cal W}_b$, to leaves of $\cal W$. This diffeomorphism respects orthogonality of two nets and preserves the coordinate foliations of bipolar coordinate system. Therefore it is conformal and, on assumption of preserving the orientations, is of the form (\ref{preservebipolar}). Let us draw  the circle of the hyperbolic pencil through $p_0$ and keep track of how the diffeomorphism $\varphi$ acts on the points of this circle. First observe that $\varphi$ can be analytically extended on the line covering the circle since the singularities of $\cal W$ are isolated. This extension is periodic due to (\ref{preservebipolar}). We claim that its period is $2\pi$. Let us parametrize the covering by the angle formed by the following two rays. The first ray is the positive ray of $x$-axis of the Cartesian coordinate system, and the second one originates from the circle center and pass through the point of the circle. As we go along the circle, the oriented tangent to hyperbolas rotates with respect to the oriented tangent to the circle in one direction on the half of the circle and then in the opposite direction on the rest of the circle. As $\varphi$ respects angles, we conclude that the period is $2\pi$ and that $\alpha=1$.
 \hfill $\Box$\\

\section{Confocal conics via linear net}
\noindent {\bf Proof of Theorem \ref{TangentTH}:} Let the forms (\ref{forms}) determine the
foliations of $\cal N$. Then the line directions of ${\cal L}$ annihilate the forms
$$
u_1= \omega_1+P\omega_2, \ \ \ \ u_2=\omega_1-P\omega_2,
$$
for some $P(x,y)$. Since the integral curves of ODEs $u_1=0$ and $u_2=0$ are rectilinear holds
$$
d\left(\frac{P-T}{1+PT}\right)\wedge u_1=0, \ \ \ \ d\left(\frac{P+T}{1-PT}\right)\wedge u_2=0,\ \ \ {\rm hence}
$$
\begin{equation}\label{PxPy}
\begin{array}{l}
P_x=\frac{(P^2+1)}{P(T^2+1)^2}[T(P^2+1)T_x+(P^2-T^2)T_y],\\
\\
P_y=\frac{(P^2+1)}{P(T^2+1)^2}[(1-P^2T^2)T_x-T(P^2+1)T_y].
\end{array}
\end{equation}
The Blaschke curvatures of two 3-webs composed of $\cal L$ and one of the foliations of   $\cal N$ are the same, their vanishing reads as
\begin{equation}\label{curvatureB}
\begin{array}{l}
\frac{T^2P^2+2P^2+1}{2T(P^2+1)}T_{xx}+\frac{2T^2P^2+P^2+T^2}{2T(P^2+1)}T_{yy}-T_{xy}=\\
\\ \ \ = \frac{T(T^2P^2+4P^2+3)T_x^2-2(2T^2-1)(P^2+1)T_xT_y+T(2T^2P^2-P^2+T^2-2)T_y^2}{T(T^2+1)(P^2+1)}.
\end{array}
\end{equation}
Modulo the equation $d(P_xdx+P_ydy)=0$, where $P_x,P_y$ are as in (\ref{PxPy}),
equations (\ref{curvatureB}) and (\ref{laplasT}) are equivalent, hence the equivalence of conditions 1. and 2. of the theorem.

The equations $d(P_xdx+P_ydy)=0$ and (\ref{curvatureB}) allow one to find $T_{xx}$ and $T_{yy}$:
\begin{equation}
\begin{array}{l}
T_{xx}=\frac{2T}{1-T^2}T_{xy}+\frac{2T(T^2-3)}{T^4-1}T_x^2+\frac{4T}{T^4-1}T_y^2+\frac{4(2T^2-1)}{T^4-1}T_xT_y,\\
\\
T_{yy}=\frac{2T}{T^2-1}T_{xy}+\frac{4T}{T^4-1}T_x^2+\frac{2T(T^2-3)}{T^4-1}T_y^2-\frac{4(2T^2-1)}{T^4-1}T_xT_y.
\end{array}
\end{equation}
Differentiating the above expressions one gets $T_{xxy}$ and $T_{xyy}$. The obtained differential   system for $P,T,T_x,T_y,T_{xy}$ turned out to be integrable.

Lines of each family of $\cal L$ either envelope a curve arc or these arcs simultaneously degenerate into points. The arcs can be parametrized (for example, by $y$) as follows:
$$
\left(x-\frac{1}{W^i_y},y-\frac{W^i}{W^i_y}\right), \ \ \  W^1=\frac{P+T}{PT-1}, \ \ \ W^2=\frac{P-T}{1+PT}.
$$

In the latter case of degeneration, these formulas give the coordinates of the vertices of the line pencils forming ${\cal L}$.  Then it is immediate from the optic properties of ellipse that all curves of one family of $\cal N$ are confocal ellipses.

In the former case, using this parametrization and the above integrable system, one checks that  both arcs belong to the same smooth conic $Q$. Then the tangents to it from any point $p$ "outside" $Q$ share the angle bisectors with the lines connecting p with the foci of $Q$. (See Theorem 2.2.4 in \cite{GSO-16}, p.42. In particular, if $Q$ is an ellipse, the claim follows from the Graves string construction.) In the case of a parabola, the ideal point
of the axis serves as the "second" focal point. One again concludes that  the leaves of $\cal N$ are  (arcs of) conics of the confocal family of $Q$, if $Q$ is not a circle, and the coordinate lines of the polar coordinate system sharing the center with the circle otherwise.

 Finally, the claims on the  rank and the rank curve follow from the following fact (see for instance \cite{BSST-17}): one can choose local first integrals $s_1,s_2$ of foliations by confocal conics so that the tangent lines  to one fixed conic $Q$ of the confocal family have equation $s_1+ s_2=const$ or $s_1- s_2=const$. (Recall that there are 2 tangents from any point "outside" $Q$.) Thus the locally defined map $(x,y)\mapsto (s_1,s_2)$ parallelizes the 4-web of $\cal N$, $
 \cal L$, all four foliations of the obtained web are families of parallel lines, the lines of each family intersect the line at infinity at a vertex, the vertices corresponding to $\cal N$ divide the vertices  corresponding to $\cal L$ harmonically. The degenerate cases when $Q$ is a parabola or a circle one can get as the limit cases.   \hfill $\Box$\\
\section{Concluding remarks}
The analytic description of confocal conics, presented in \cite{BSST-17}, provides more examples of "geometric" 4-webs of maximal rank.  For instance, the 4-web, composed of confocal family of conics, hyperbolic pencil of circles sharing foci with conics, and straight lines orthogonal to the line joining the foci, is also parallelizable (see subsection 8.6.1 of \cite{BSST-17}). In this paper, we have restricted ourselves with the most "symmetric" cases.

The rich web structures related to confocal conics invite questions on the rank of webs composed of even more foliations in hope to find some "master" web of maximal rank including all 4-webs considered above as subwebs. A straightforward try failed: the 6-web composed of confocal conics, Cartesian coordinate lines and bipolar coordinate lines is not of maximal rank.

\section*{Acknowledgements}
The author thanks A. Bobenko and W. Schief for useful discussions.
This research was supported by FAPESP grant \#2018/20009-6 and partially by SFB/TRR 109 "Discretization in Geometry and Dynamics". The author also thanks  the personnel of the Institute of Mathematics of Technische Universit\"at Berlin, where this study was initiated, for their warm hospitality.

\end{document}